\documentclass[10pt,a4paper]{amsart}
\usepackage[utf8]{inputenc}
\usepackage[T1]{fontenc}
\usepackage{amsmath}
\usepackage{amssymb}
\usepackage{graphicx}
\usepackage{mathrsfs}
\usepackage{hyperref}
\usepackage{cleveref}
\hypersetup{
	colorlinks,
	linkcolor={red!50!black},
	citecolor={blue!50!black},
	urlcolor={blue!80!black}
}
\usepackage{tikz-cd}
\usepackage[margin=3.5cm]{geometry}

\newtheorem{thm}{Theorem}[section]
\newtheorem{cor}[thm]{Corollary}
\newtheorem{lem}[thm]{Lemma}
\newtheorem{prop}[thm]{Proposition}
\theoremstyle{definition}
\newtheorem{rem}[thm]{Remark}
\newtheorem{definition}[thm]{Definition}
\newtheorem{ex}[thm]{Example}

\newtheorem{thmintro}{Theorem}

\newcommand{\HH}{\mathbb{H}}
\newcommand{\Z}{\mathbb{Z}}
\newcommand{\C}{\mathbb{C}}

\newcommand{\cO}{\mathcal{O}}
\newcommand{\cN}{\mathcal{N}}
\newcommand{\cD}{\mathcal{D}}

\newcommand{\R}{\mathbb{R}}
\newcommand{\CP}{\mathbb{CP}}
\newcommand{\id}{\operatorname{id}}
\newcommand{\U}{\mathrm{U}}
\newcommand{\SU}{\mathrm{SU}}
\newcommand{\Sp}{\mathrm{Sp}}
\newcommand{\SO}{\mathrm{SO}}

\newcommand{\oX}{\overline{X}}

\newcommand{\bq}{\bar{q}}

\newcommand{\fc}{\mathfrak{c}}

\newcommand{\dY}{\delta_Y}

\newcommand{\dW}{\delta_W}

\newcommand{\supp}{\operatorname{supp}}
\newcommand{\del}{\partial}
\newcommand{\delbar}{\bar\partial}
\usepackage{tikz-cd}

\title{Holomorphic linking numbers, ABC Massey products, and Calabi-Yau 3-folds}

\author{Luc\'ia Mart\'in-Merch\'an}
\address{Institut für Mathematik, Humboldt Universität zu Berlin}
\email{lucia.martin.merchan@hu-berlin.de}

\author{Jonas Stelzig}
\address{Mathematisches Institut der Ludwig-Maximilians-Universität München}
\email{jonas.stelzig@math.lmu.de}

\begin{document}

\maketitle

\begin{abstract}
\noindent On compact Kähler manifolds, we relate ABC Massey products arising from complex analytic cycles to holomorphic linking numbers.
This enables us to construct a family of simply connected projective 3-folds with trivial canonical bundle, equipped with a non-vanishing ABC Massey product.
\end{abstract}
\section{Introduction}

Compact Kähler manifolds are rationally formal \cite{DGMS}. In particular, all Massey products must vanish. This result covers in some sense 'most' of the known oriented compact Riemannian manifolds with special holonomy, namely those with restricted holonomy group contained in $\U(n)$, $\SU(n)$, and $\Sp(n)$. It had inspired the folklore conjecture that any compact manifold admitting a Riemannian metric with restricted holonomy strictly contained in $\SO(n)$ should be formal. However, despite a few positive results, e.g. \cite{AmaKap_fibr}, \cite{Crowley2020}, a counterexample to this conjecture was quite recently found in \cite{MartinMerchan2024}.

Pluripotential homotopy theory \cite{StePHT} is a refinement of rational homotopy theory on complex manifolds which takes into account holomorphic information. It gives rise to a stronger notion of formality, which is obstructed by the non-vanishing of a secondary version of Massey products, called ABC Massey products \cite{AnToBCform}, \cite{MilSt_bigrform}. On a compact Kähler manifold $X$, these operations can be thought of as partially defined operations
\[
\begin{tikzcd}
\langle \_,\_,\_\rangle_{ABC}:H_{dR}(X)^{\otimes 3}\ar[r, dashed]&  H_{dR}(X)
\end{tikzcd}
\]
of degree $(-1,-1)$ with respect to the bigrading given by the Hodge decomposition. These products are also very closely related to Deninger's products in Deligne cohomology \cite{Den_higher}.

As shown in \cite{PSZ24}, in contrast to \cite{DGMS}, the ABC Massey products can be nontrivial quite often even on compact Kähler (or even  projective) manifolds. However, the examples in \cite{PSZ24} do not have restricted holonomy contained in $\SU(n)$ or $\Sp(n)$, and so it was left open whether Calabi-Yau or Hyperkähler manifolds are strongly formal.

In this paper we answer the Calabi-Yau case in the negative as follows:

\begin{thmintro}\label{thmintro: ABC MP on CY}
    There is a family $\{X_{\tau}\}_{\tau\in \HH}$ over the upper half-plane of simply connected, projective $3$-folds with trivial canonical bundle and on each $X_\tau$ there exist divisors $D_1,...,D_4$ on $X_\tau$ such that
\[
\int_{D_1}\langle D_2,D_3, D_4\rangle_{ABC}=\frac{4}{\pi}\cdot \log|\lambda(\tau^{-1})|~,
\]
where $\lambda$ is the modular $\lambda$-function. In particular, for general $\tau$, one has $\langle D_2, D_3, D_4\rangle_{ABC}\neq 0$ and the manifold $X_\tau$ is not strongly formal.
\end{thmintro}
Our proof combines techniques and ideas from \cite{MartinMerchan2024} and \cite{PSZ24}. As a byproduct which may be of independent interest, we relate ABC Massey products to holomorphic linking numbers. Recall that the holomorphic linking number of two disjoint, null-homologous cycles $Z,W$ on a compact Kähler manifold is given by

\[
\langle Z,W\rangle =\int_Z Q_W~,
\]
where $Q_W$ is a current, smooth outside the support of $W$, s.t. $i\del\delbar Q_W=\delta_W$, for the current of integration $\delta_W(\omega)=\int_W\omega|_W$. This definition is analogous to a formula for the classical linking number of two null-homologous cycles. It also arises in Arakelov geometry as the Archimedean component of the height pairing, see e.g. \cite{Beilinson_HeightPairing}, \cite{Bloch_HeightPairing}, \cite{Bost_HeightPairing}, \cite{HainBiext}. Moreover, it is computable by geometric means. We show:

\begin{thmintro}\label{thmintro: MP vs LN}
   Let $A, B, C, D \in Cyc(Y)$ be analytic cycles of codimensions $a, b, c, d$, respectively, with $a + b + c + d = n + 1$. 
   Let $F,F',G,G'\in Cyc(Y)$  be null-homologous of codimensions $f,f',g, g'$ with $f+g'=f'+ g=n+1$ such that
    \begin{enumerate}
        \item $A\cap  B\cap C=\emptyset, A \cap C \cap D=\emptyset$.
        \item $|F|\subseteq |A|\cap |B|$, $|F'|\subseteq |A|\cap |D|$,
        $|G|\subset |B|\cap |C|$, $|G'|\subset |D|\cap |C|$.
        \item The components of 
        $A=\sum n_j A_j$ and $C=\sum m_k C_k$ are disjoint and smooth, and the following equalities hold in $H^*(|A|)$ and $ H^*(|C|)$ respectively:
  \begin{align*}
        &[A]_{rel}\cdot [B]|_A=[F]_{rel},  &[A]_{rel}\cdot [D]|_{A}=[F']_{rel},  \\
        &[C]_{rel}\cdot [B]|_{C}=[G]_{rel},  &[C]_{rel}\cdot [D]|_{C}=[G']_{rel}. 
        \end{align*}
    \end{enumerate}
    Then, there is an equality:
    \[
        \int_A\langle B,C,D\rangle_{ABC}=\langle F,G'\rangle- \langle F',G\rangle.
    \]
\end{thmintro}

In this theorem, we denoted $[F]_{rel}\in H(|A|)$ for the class whose restriction to each component $A_j$ is $\sum_{F_k\subset A_j} m_k [F_k]$.

\noindent\textbf{Related work:} In the forthcoming PhD thesis of A.~P. Contr\`o, it will be shown that Hyperk\"ahler manifolds are indeed strongly formal.

\section{Holomorphic linking numbers and ABC Massey products}
Throughout this section, $Y,Y'$ will be  connected complex manifolds of dimensions $n,n'$ (reserving the letter $X=X_\tau$ for the examples in Theorem \ref{thmintro: ABC MP on CY}). Our main interest is in the case $Y$ compact and K\"ahler.

\subsection{Preliminaries} 
We recall some notions and introduce some notations. We mainly follow \cite{Bost_HeightPairing}, \cite[section 1]{GilletSoule}, \cite{Hoermander_LPDO_I}. 

\subsubsection{Forms and currents}We write $A^{p,q}(Y)$, resp. $A_c^{p,q}(Y)$, for the space of $\C$-valued smooth, resp. smooth compactly supported, differential forms of type $(p,q)$ and $A(Y)=\bigoplus_{p,q} A^{p,q}(Y)$, resp. $A_c(Y)=\bigoplus A^{p,q}_c(Y)$ for the bicomplexes of all (compactly supported) forms, with differential $d=\del+\delbar$. When $Y$ is compact, $A(Y)=A_c(Y)$. 
Similarly, we write $\cD^{p,q}(Y)$ (resp. $\cD_c^{p,q}(Y)$) for the space of (compactly supported) currents of degree $(p,q)$, i.e. the topological dual of 
$A^{n-p,n-q}_c(Y)$ (resp. $A^{n-p,n-q}(Y)$), and $\cD(Y)=\bigoplus\cD^{p,q}(Y)$ (resp. $\cD_c(Y)=\bigoplus\cD^{p,q}_c(Y)$) for the spaces of all currents. Also $\cD (Y)$ and $\cD_c(Y)$ are bicomplexes, equipped with differential $T\mapsto ( \omega\mapsto (-1)^{|\omega|+1}T(d \omega))$. 
The inclusion
$$
A(Y) \hookrightarrow \cD(Y), \qquad \alpha \mapsto (\delta_Y \wedge \alpha) (\beta):=\int_{Y}{\beta \wedge \alpha}
$$
is a map of bicomplexes by Stokes' theorem and a map of $A(Y)$-modules, where one defines the module structure on the right by  $T\wedge \alpha (\beta)= T(\beta \wedge \alpha)$ for a for a given current $T$. This inclusion is a pluripotential quasi-isomorphism: It induces isomorphisms in Dolbeault, Bott-Chern, Aeppli and de Rham cohomology.


The expression 
$\delta_Y \wedge \eta$ still defines a current when $\eta\in L^1_{\mathrm{loc}}(Y)$, i.e. when the coefficients of $\eta$ in local holomorphic coordinates are measurable functions defined almost everywhere, whose modulus is locally integrable.

A proper holomorphic map $f\colon Y\to Y'$,  induces a pullback map of bicomplexes $f^*:A_c(Y')\to A_c(Y)$ and a pushforward
$f_*\colon \cD(Y) \to \cD(Y')[r]$, where $r=n'-n$ and  $(\cD(Y')[r])^{p,q}=\cD^{p-r,q-r}(Y')$, via
$f_*(T)(\beta)= T(g^*\beta)$.

Moreover, if $f\colon Y\to Y'$ is a surjective holomorphic submersion, then integration along the fibers defines a bicomplex map $f_* \colon A_c(Y)\to A_c(Y')[r]$.
Taking the dual gives pullback maps $f^*:\cD(Y)\to\cD(Y')$. When $f$ is both proper and a submersion, the operations on forms and currents are compatible in that $f_*(\delta_Y \wedge \alpha)=\delta_{Y'}\wedge f_*(\alpha)$ and $f^*(\delta_Y'\wedge \alpha)=\delta_Y\wedge f^*\alpha$. 

Another instance where pullback maps on currents are defined is the inclusion of open subsets $i:U\subseteq Y$. The inclusion $i_*:A_c(U)\to A_c(Y)$ induces $i^*:\cD(Y)\to \cD(U)$, which we denote by $T\mapsto i^*T=:T|_U$. We say $T$ is smooth on $U$ if there exists $\alpha\in A^{p,q}(U)$ with $T|_U=\delta_U\wedge \alpha$. The support $\supp(T)$ is the set of all points $y\in Y$ such that $T|_V\neq 0$ for any open neighborhood $y\in V\subseteq Y$. The singular support $\mathrm{sing}\supp(T)$ is the set of all $x\in Y$ for which there is no open neighborhood $x\in V$, the restriction $T|_V$ is smooth.

\subsubsection{Analytic cycles} An (analytic) cycle is a locally finite formal sum $Z=\sum_{i\in I} n_i Z_i$ where  $\{ Z_i\}_{i\in I}$ are different $d$-dimensional irreducible analytic subsets, called the components of $Z$, and $n_i\in \Z$, called multiplicities. We denote by $Cyc_k(Y)$, resp. $Cyc^ k(Y)$, the abelian group of analytic cycles of pure dimension, resp. codimension, $k$ on $Y$. Given a cycle $Z=\sum n_i Z_i$ we denote by $|Z|=\bigcup Z_i$ the support, i.e. the underlying topological space of $Z$. There are proper pushforward and flat pullback operations of currents: Let $f:Y\to Y'$ be a holomorphic map.

\begin{enumerate}
    \item If $f$ is proper, there is a group homomorphism $f_*:Cyc_k(Y)\to Cyc_k(Y')$ given on an irreducible analytic subset $Z\subseteq X$ of dimension $k$ as
    \[
        f_*(Z)=\begin{cases}
            \deg(f|_Z)\cdot f(Z)&\text{if }f|_Z:Z\to f(Z)\text{ is generically finite}\\
            0&\text{else.}
        \end{cases}
    \]
    and extended additively. The degree is the number of preimages in $Z$ over a generic point of $f(Z)$, or equivalently, the degree of the extension of function fields $k(f(Z))\subseteq k(Z)$.
    
    \item If $f$ is flat (open), there is a pullback map $f^*:Cyc^k(Y')\to Cyc^{k}(Y)$ given on an irreducible analytic subset by $Z'\subseteq Y'$ of dimension $k$ by
     \[
        Z'\mapsto f^*Z':=\sum_{\substack{Z\subseteq f^{-1}Z'\\\text{irreducible components}}}m_{f^{-1}Z', Z}\cdot Z
     \]
     where $m_{f^{-1}Z',Z}:=\mathrm{length}_{\cO_{f^{-1}Z',Z}}\cO_{f^{-1}Z',Z}$ equals $1$ if $f$ is an open inclusion or, more generally, a submersion near $f^{-1}Z$.
     \item If $f$ is an embedding, there is a pullback map for analytic codimension-$k$ cycles whose irreducible components $Z'$ intersect $Y$ properly; i.e,  $\dim(f(Y)\cap Z)= n-k$. For such an irreducible component $Z'$, 
     $$  f^*(Z'):=\sum_{\substack{Z\subseteq f^{-1}Z'\\\text{irreducible components}}}Z~ \in Cyc^{k}(Y).
     $$
\end{enumerate}

\subsubsection{Cycles and currents}
We now assume $Y,Y'$ to be compact K\"ahler manifolds of complex dimensions $n,n'$.
Given a cycle $Z=\sum_{i=1}^ N n_i Z_i\in Cyc^k(Y)$, we denote by $\delta_Z=\sum_{i=1}^N n_i \delta_{Z_i}$ its Dirac-delta current, where for any component $\delta_{Z_i}$ is given by $ \alpha\mapsto \int_{Z_i}\alpha:=\int_{Z_i^{reg}}\alpha$. The current $\delta_Z$ is $\del$ and $\delbar$-closed and we denote the cohomology class by $[Z]:=[\delta_Z]\in H^{k,k}(Y)$ and call it the fundamental class of $Z$. It is the Poincar\'e dual to a homological fundamental class for $Z$.

We say that $T\in \cD^{k-1,k-1}(Y)$ is a Green current for an analytic cycle $Z\in Cyc^k(Y)$ if $\del\delbar T- \delta_Z$ is smooth on $Y$. Since the inclusion $A(Y)\to \cD(Y)$ is a pluripotential quasi-isomorphism, Green currents always exist.

The Lelong-Poincaré formula allows us to find Green currents for divisors. Given a divisor $Z$, we consider $L=\cO(Z)$ the line bundle associated to $Z$ and a hermitian metric $h$ on $L$. Let $s$ be a holomorphic section with $\mathrm{div}(s)=Z$, and let $\omega$ be the curvature of the Chern connection of $(L,h)$, which represents $c_1(L)=[Z]$. Explicitly, if on a trivializing open set $U$ we have $h|_{L|U} = e^{-\varphi}|\cdot|$, then $\omega|_{L|_U}=\frac{i}{2\pi}\del \delbar \varphi$. The Lelong-Poincarè formula shows that $\frac{i}{\pi}\log \|s\|_h$ is a Green current for Z because
\begin{equation} \label{eqn:Green-currents-for divisors}
\frac{i}{\pi}\del \delbar (\delta_Y \wedge\log \|s\|_h)=\delta_Z - \omega.
\end{equation}
In fact, we can assume that $\|s\|_h=1$ outside an open neighborhood $U$ of $Z$, so that both $\omega$ and $ \log(\|s\|_h)$ are supported on $U$. 

Using the projection formula, one sees that pullback and pushforward operations on currents and cycles are compatible in the sense that $f^*\delta_Z=\delta_{f^*{Z}}$ for $f$ a surjective holomorphic submersion or an open inclusion and $f_*\delta_Z=\delta_{f_*Z}$ for $f$ a proper holomorphic map.

\subsubsection{Wave front sets} For every current $T\in \cD^k(Y)$, the wave-front $WF(T)$ is a subset of $WF(T)$ of $T^\ast Y- Y$, measuring the directions in which $T$ is singular. To define it, we first recall that the Fourier transform of a temperate distribution $u$ (i.e, in the topological dual of Schwartz functions)  is given by 
$$
\widehat{u}(f)= u(\widehat{f}), \, f\in C^{\infty}_c(\R^m), \, \widehat{f}(\xi)=\int_{\R^m} e^{-i \langle x,\xi\rangle}f(x)dx.
$$
 If $u \in \cD^m_c(\R^m)$ then  
 $\widehat{u} \in C^{\infty}(\R^m)$ and $\widehat{u}(\xi)=u(\mathrm{e}_\xi)$ where $\mathrm{e}_{\xi}(x)= e^{-i\langle x,\xi\rangle}$ (the left-hand side is well-defined because $u$ has compact support; see \cite[Theorem 7.1.14]{Hoermander_LPDO_I} for details). In this case, we say that $\widehat{u}$ has a rapid decay along $\xi_0$ if there is a neighborhood $V$ of $\xi_0$, conic with respect to the multiplication by positive scalars, such that for every $d\in \mathbb{N}$ and $\xi \in V$ one has
$$\sup_{\xi \in V}(1+\|\xi\|)^d|\widehat{u}(\xi)|<\infty.
$$

Let $U\subset \R^m$ be an open set, $u \in \cD^m(U)$ and $(x_0,\xi_0)\in U\times \R^m-\{0\}$ then
$(x_0,\xi_0)\notin WF(u)$  if and only if there is $\phi\in C^{\infty}_c(U)$ with $\phi(x_0)\neq 0$ such that $\widehat{\phi u}\in C^{\infty}(\R^m)$ has rapid decay at $\xi_0$.
This definition is first extended to distributions on manifolds (see \cite[p. 256]{Hoermander_LPDO_I}), and then to currents, as these are locally forms with distributions as coefficients.
We will need some properties, which generalize those for distributions and can be found in \cite{Hoermander_LPDO_I, ReedSimon}:
\begin{enumerate}
    \item If $p:T^*Y\to Y$ is the projection, $p(WF(T))=\mathrm{sing}\supp T$ \cite[Theorem IX.44 (d)]{ReedSimon}.
     \item If $S,T\in\cD(Y)$ and $f\in C^\infty(Y)$, then $WF(f T)\subset WF(T)$ and $WF(S+T)\subseteq WF(S)\cup WF(T)$ \cite[equation (8.1.9), p.256]{Hoermander_LPDO_I}, \cite[Theorem IX.44 (c)]{ReedSimon}.
    \item If $Z\subseteq Y$ is a submanifold $WF(\delta_Z)$ is the conormal bundle of $Z\subseteq Y$(without the zero section).
    \item If $T$ is a real current, i.e $T(\bar{\alpha})=\overline{T(\alpha)}$, then
     $WF(T)=-WF(T)$. This property follows from
    $\widehat{\bar{\phi} u}(-\xi)=\overline{ \widehat{\phi u}(\xi)}$ if $u\in \cD^m(U)$ is a real distribution and $\phi \in C^\infty_c(U)$. In particular, if $Z\subset Y$ is an analytic cycle then $WF(\delta_Z)=-WF(\delta_Z)$
    \item If $P$ is a linear partial differential operator with smooth coefficients and $T\in \cD(Y)$ a current, then $WF(PT)\subseteq WF(T)$ \cite[equation (8.1.11),  p. 256]{Hoermander_LPDO_I}. For $P$ elliptic, $WF(PT)=WF(T)$ \cite[Corollary (8.3.2)]{Hoermander_LPDO_I}.
    \item If $S,T\in\cD(Y)$ and  $WF(T)\cap -WF(S)=\emptyset$ , then the product $S\wedge T$ is defined and $WF(S\wedge T)\subseteq WF(S)\cup WF(T) \cup W$, where $W=\{ \xi_y+\xi'_y : \xi_y\in WF(T), \, \xi_y'\in WF(S) \}$ \cite[Theorem 8.2.10]{Hoermander_LPDO_I}. It extends the product of forms with currents and satisfies the Leibniz rule.
    \item Let $f \colon Y  \to Y'$  be a holomorphic map, and define the conormal set as $N_f=\{ \alpha_y \mbox{ s.t. } y \in Y, 0\neq \alpha_y \in T_{f(y)}^*Y, \alpha_y \circ df_y=0\}$. For any current $T\in \cD^{k}(Y')$ there is a unique pullback $f^*T\in \cD(Y)$ that continuously extends the usual pullback of currents determined by smooth forms, satisfying  $WF(f^*T)\subset f^*WF(T)$, and commuting with $\del,\delbar$  and with the product of currents \cite[Theorem 8.2.4]{Hoermander_LPDO_I}.
\end{enumerate}

We also refer to \cite{Bost_HeightPairing} for the rephrasing (without proofs) in the global setting.

\subsection{Holomorphic linking numbers}
In this section, $Y.Y'$ are always compact connected K\"ahler manifolds of dimensions $n,n'$. We recall the definition of holomorphic linking numbers and prove two properties. These may be known to experts, but we could not find a reference.

\begin{definition}[E.g. \cite{Beilinson_HeightPairing}, \cite{HainBiext}]
    The holomorphic linking number of two disjoint null-homologous analytic cycles $Z\in Cyc^d(Y)$, $W\in Cyc^e(Y)$, s.t. $d+e=n+1$ is given by 
\[
\langle Z,W\rangle:=\int_{Z}\eta_W=\int_W \eta_Z~,
\]
where $\eta_Z ,\eta_W$ are smoth forms on $Y-|Z|,Y-|W|$ of pure bidegree, and $T_Z,T_W$ are currents such that $i\del\delbar T_Z=\delta_Z$, $i\del\delbar T_W=\delta_W$, and  $T_Z|_{Y-|Z|}= \delta_{Y-|Z|}\wedge \eta_Z$, $T_W|_{Y-|W|}=  \delta_{Y-|W|}\wedge \eta_W$.
\end{definition} 
In fact, our definition differs by a factor of $\pi$ from that in \cite{HainBiext}. In addition, since $WF(T_W)\cap -WF(\delta_Z)=\emptyset$, we can rewrite $\langle Z,W\rangle = (\delta_Z \wedge T_W)(1)$; the equality $(\delta_W \wedge T_Z)(1)=(\delta_Z \wedge T_W)(1)$ follows from the Leibniz rule applied to $T_Z \wedge T_W$.

\begin{ex} \label{ex:linking-Riemann-Surface}
Let $\Sigma$ be a compact connected Riemann surface and $f,g\colon \Sigma \to \CP^1$ be meromorphic functions, and let $\mathrm{div}(f),\mathrm{div}(g)$ be their divisors of zeros and poles. Assume that $|\mathrm{div}(f)|\cap |\mathrm{div}(g)|=\emptyset$, and write  $\mathrm{div}(g)=\sum_{k=1}^m a_k P_k$.

The holomorphic line bundle $\mathcal{O}(\mathrm{div}(f))$ is trivial. Taking $s=f$ and $h=|\cdot |$  in the Lelong-Poincaré formula \eqref{eqn:Green-currents-for divisors}, we obtain that $\frac{1}{\pi}  (\delta_Y \wedge \log |f|)$ is a $i\del \delbar$-primitive of $\delta_{\mathrm{div}(f)}$. Hence,
$$ \langle \mathrm{div}(f), \mathrm{div}(g) \rangle=\frac{1}{\pi} \log \left|f(P_1)^{a_1}\cdots f(P_m)^{a_m}\right|. $$

Applying this to $P,Q,R,S \in \Sigma=\CP^1$ and $f,g\colon \CP^1\to \CP^1$, 
\[f(z)= \dfrac{z-P}{z-Q},\quad g(z)=\dfrac{z-R}{z-S},\]
and denoting the cross ratio of the given points by  $(P,Q,R,S)$, we obtain
$$
\langle P-Q,R-S\rangle = \frac{1}{\pi}\log \left| \dfrac{(R-P)(S-Q)}{(R-Q)(S-P)} \right|= \frac{1}{\pi} \log|(P,Q,R,S)|.
$$ 
\end{ex}

Lemma \ref{lem:flexibility} below shows there is more flexibility in choice of the primitive $T$ in the preceding definition. As a consequence, the holomorphic linking number is very computable (see also \Cref{sec: defn and evaluation}). It also recovers that $\langle Z,W\rangle$ is well-defined. To prove it, we first need the following result, which generalizes \cite[Proposition 1.1]{Bost_HeightPairing}.
\begin{prop} \label{prop:regularity}
Let $T \in \cD(Y)$ with $\del \delbar T=0$. There is a harmonic form, $\omega$ and $T_1,T_2\in \cD(Y)$ with $WF(T_1),WF(T_2)\subset WF(T)$ such that 
$$ T= \delta_Y \wedge \omega + \del T_1 + \delbar T_2. 
$$
In addition, if $S\in \cD(Y)$, satisfies  $\del\delbar S=0$ and $[T]=[S]$, there is $T_3\in \cD(Y)$ with $WF(T_3)\subset WF(T)\cup WF(S)$ and $T-S=\del\delbar T_3$
\end{prop}
\begin{proof}
Consider the Laplacian $\Delta:=\Delta_\del=\Delta_{\delbar}=\frac{1}{2}\Delta_d$ and $G$ be the Green's operator for $\Delta$. In the sequel, we use that both $\Delta,G$ commute with $\del,\delbar,\del^*,\delbar^*$, and that $H=\id - G \Delta$ is the orthogonal projection to the spaces of harmonic forms (\cite[Sections 0.6,0.7]{GH}). 
These operators extend naturally to currents, and similar properties hold;  in particular, since $\Delta H=0$, by elliptic regularity of the Laplacian, $HT$ is smooth for any $T\in \cD(Y)$.

 Given a current $T$, we have
$T= HT + G\Delta T = HT + \delbar (\delbar^*GT) +\delbar^*(\delbar GT)$. In addition,
$$  
\delbar^* \delbar G T= H(\delbar^* \delbar GT) + \Delta G^2 \delbar^* \delbar T=  H(\delbar^* \delbar GT) + \del (\del^*  \delbar^* G^2 \delbar T).
$$
In the last step we used $\del \delbar T=0$ and that $\del$ commutes with $\delbar^*$. Since $H$ is a smoothing operator, there is a harmonic form $\omega$ with
$ H(T +\delbar^* \delbar GT )=\delta_Y \wedge \omega$; in addition, the computations above imply that we can take 
$$ T_1=\del^*  \delbar^* \delbar G^2 T,\qquad  T_2=\delbar^*GT.
$$
Using the first part of property (3) for wave front sets: $WF(T_1)\subset WF(G ^2 T)$, $WF(T_2)\subset WF(GT)$. The proof follows from the property $WF(G T')\subset WF(T')$ which we check using that $HT$ is smooth, the operator $\Delta$ is elliptic, and  properties 
(3),(4):
$WF(GT')=WF(\Delta G T') = WF(T'-HT')\subset WF(T')$.

For the second part, we note that since $[S]=[T]$, we have $T-S=\del (T_1-S_1)+ \delbar(T_2-S_2)$ with $WF(S_1),WF(S_2)\subset WF(S)$. Using the  Hodge decomposition for these, we obtain for $j=1,2$:  $T_j-S_j= \delta_Y \wedge \omega_j + \del T_3^j + \delbar T_3^j$ with $WF(T_3^j)\subset WF(T)\cup WF(S)$. Hence $T-S= \del \delbar T_3 $ for $T_3=T_3^1- T_3^2$ and $WF(T_3)\subset WF(T)\cup WF(S)$.
\end{proof}

\begin{lem} \label{lem:flexibility}
Let $Z\in Cyc^d(Y)$, $W\in Cyc^e(Y)$ be nullhomologous with $d+e=n+1$.\\
For any $T\in \cD^{d-1,d-1}(Y)$ with $WF(T)\cap WF(\delta_W)=\emptyset$ and $i \del \delbar T= \delta_Z$ there is an equality
\[\langle Z,W\rangle=(\dW\wedge T)(1).\]
\end{lem}
\begin{proof}
Let $T_Z\in \cD^{d-1,d-1}(Y)$, smooth on $Y-|Z|$ such that $i \del \delbar T_Z = \delta_{Z}$ and $WF(T_Z)\subset WF(\delta_Z)$.
By Proposition \ref{prop:regularity}, 
$
T-T_Z= \delta_Y\wedge \omega + \del T_1+ \delbar T_2
$
with $\omega\in A(Y)$ harmonic and $WF(T_1),WF(T_2)\subset WF(T-T_Z)\subset WF(\delta_Z)\cup WF(T)$. In particular, $WF(T_j)\cap -WF(\delta_W)=\emptyset$ for $j=1,2$, hence the products $\delta_W\wedge T_j$ are well-defined. Using the Leibniz rule and $\del \delta_W=\delbar \delta_W=0$, we obtain:
$$
\left( \dW \wedge \del T_1 + \dW \wedge \delbar T_2 \right)(1)= \del( \dW \wedge T_1 ) (1) + \delbar(\dW \wedge T_2)(1)=0. 
$$
Since $\omega$ is harmonic on $Y$, and $W$ is null-homologous, we also have $ (\delta_W \wedge (\delta_Y \wedge \omega))(1)=\int_W \omega=0$. Hence,
$ \langle Z, W \rangle - (\dW \wedge T) (1) =   (\dW \wedge (T-T_Z)) (1)= 0$.
\end{proof}

The following is a functoriality property for the holomorphic linking number. 

\begin{prop} \label{prop:adjunction}
Let $p:Y\to Y'$ be a holomorphic map, and let 
$Z\in Cyc^d(Y)$, $W\in Cyc^{e}(Y')$ be null-homologous with $d+e=n+1$ with $p(Z)\cap W=\emptyset$. Assume one of the following conditions:
\begin{enumerate}
\item The map $p:Y\to Y'$ is surjective, $\dim(p(Z))=\dim(Z)$, and there is an open neighborhood $|W|\subseteq U'\subseteq Y'$ with $p(|Z|) \cap U' = \emptyset$ such that $p|_{p^{-1}U'}$ is a submersion. 
\item The map $p\colon Y \hookrightarrow Y'$ is a submanifold embedding and the components of $W$ are smooth and intersect $f(Y)$ transversally.
\end{enumerate}
Then, 
\[\langle Z, p^{*}(W)\rangle =  \langle p_*(Z),W \rangle.\]
\end{prop}
\begin{proof}
We first note that indeed $[p(Z)]=p_*[Z]=0$; to prove $[p^{*}(W)]=0$, we let 
 $\delta_W=i\del\delbar T$. We can assume that $WF(T)\subset WF(\delta_W)$ by Proposition \ref{prop:regularity}, taking $S=0$.
 We first prove that in both situations, $WF(\delta_W)\cap N_p=\emptyset$, which implies that $p^*T$, $p^*\delta_W$ are well-defined. In the first case, set $U=p^{-1}(U')$, then $WF(\delta_W)\subset T^*Y'|_{U'}$ and 
 $N_p \cap  T^*Y'|_{U'}=\emptyset$ because $p|_{U}$ is a submersion. In the second, if $W=\sum n_j W_j$, then $WF(\delta_W)= \cup_{j} (N^*W_j-W_j)$. Since $N_p=N^*p(Y)-p(Y)$ and the intersection of $W_j$ and $p(Y)$ is transverse, $(N^*W_j -W_j)\cap (N^*p(Y)-p(Y))=\emptyset$.

 Then $i \del \delbar p^*T=p^*\delta_W$. We now verify that $p^*\delta_W = \delta_{p^{-1}W}$. In the first case, this follows from the uniqueness of the pullback extension and the fact that, for a submersion, the pullback corresponds to integration along the fibers. More precisely, $p^*\delta_W(\beta)=p^*\delta_W(\beta|_U)=\delta_W(p_*(\beta|_U))= \delta_{p^{*}W} (\beta).$  In the second, \cite[Example 8.2.8]{Hoermander_LPDO_I} shows that, for every irreducible component $W_j$  we obtain $p^*\delta_{W_j}= \delta_{p^{-1}(W_j\cap Y)}$; and the statement follows from this.

 Hence $i\del\delbar p^*T=\delta_{p^*W}$ and $[p^*W]=[0]$. In addition $WF(p^*T)\subset p^*WF (\delta_W)$ and $p^*WF(\delta_W) \cap WF(\delta_Z)=\emptyset$ because $p^{-1}(|W|)\cap |Z|=\emptyset$; indeed, $T|_{Y-p^{-1}(|W|)}=\delta_{Y-p^{-1}(|W|)} \wedge \eta$ for some smooth form $\eta$. Hence,  
 \[\langle Z,p^*W \rangle = \int_Z p^*\eta = \int_{p_*(Z)} \eta = \langle p_*(Z), W \rangle. 
 \]
\end{proof}

\begin{rem}
We have only used the special case of \Cref{lem:flexibility} where $T$ is nonsingular on some neighborhood of $W$. Using more one can likely weaken the conditions on $p$ here, but we will not need this.
\end{rem}

\subsection{Linking numbers and Massey products} $Y$ continues to be a connected, compact K\"ahler manifold of dimension $n$. Our next goal is to connect (certain) ABC Massey products to holomorphic linking numbers. This is analogous to the topological case, cf. \cite[Prop. 3]{MartinMerchan2025}. Let us first recall the definition of ABC Massey products from \cite{AnToBCform}:

\begin{definition}
Let $\fc_1, \fc_2,\fc_3\in H_{BC}(Y)$ with $\fc_1\wedge \fc_2=\fc_2\wedge\fc_3=0\in  H_{BC}(X)$.
\begin{enumerate}
    \item A defining system is a $5$-tuple of smooth forms $(\alpha,\beta,\gamma,x,y)$ where $\alpha,\beta,\gamma$ are $\del$- and $\delbar$-closed such that \[
    \fc_1=[\alpha],~\fc_2=[\beta],~\fc_3=[\gamma]
    \] and 
    \[i\del\delbar x=\alpha\wedge \beta,\quad \del\delbar y = \beta\wedge\gamma.
    \]
    \item The ABC Massey product is the coset
    \[
\langle\fc_1,\fc_2,\fc_3\rangle_{ABC}:=[\alpha y-x\gamma]\in H_A(Y)/([\alpha]H_A(Y)+H_A(Y)[\gamma]),
\]
where $(\alpha,\beta,\gamma, x,y)$ is a defining system.
\end{enumerate} 
\end{definition}
\begin{rem}\,
\begin{enumerate}
    \item As notation suggests, the ABC Massey product $\langle \fc_1,\fc_2,\fc_3\rangle_{ABC}$ depends only on the classes $\fc_1,\fc_2,\fc_3$ and not on the choice of a defining system.
    \item If  $D,E,F\in Cyc(Y)$ s.t. $[D]\wedge [E]=[E]\wedge [F]=0\in H_{dR}(Y)$ we will abuse notation and write $\langle D,E,F\rangle_{ABC}:=\langle [D],[E],[F]\rangle_{ABC}$.
    \item The definition is phrased in terms of Bott-Chern and Aeppli cohomology since in that way it is meaningful on any complex manifold. Since we assume $Y$ to be compact K\"ahler, the natural maps \[
    H_{BC}(Y)\to H_{dR}(Y)\to H_{A}(Y)
    \]
    induced by the identity are isomorphisms and one could work with de Rham cohomology throughout. 
\end{enumerate}
\end{rem}

If $A=\sum n_j A_j\in Cyc(Y)$ is a cycle with smooth disjoint components and $F=\sum m_k F_k$ another cycle with $|F|\subseteq |A|$, then $F$ defines a class $[F]_{rel}\in H(|A|)=\bigoplus H(A_j)$ via $[F]_{rel}|_{A_j}=\sum_{F_k \subset A_j} m_k [F_k]$. Denoting the inclusion map $j:|A|\hookrightarrow Y$, we have $j_*[F]_{rel}=[F]$. 

\begin{thm}[= \Cref{thmintro: MP vs LN}]\label{thm: MP vs LN}
   Let $A, B, C, D \in Cyc(Y)$ be analytic cycles of codimensions $a, b, c, d$, respectively, with $a + b + c + d = n + 1$. 
   Let $F,F',G,G'\in Cyc(Y)$  be null-homologous of codimensions $f,f',g, g'$ with $f+g'=f'+ g=n+1$ such that
    \begin{enumerate}
        \item $A\cap  B\cap C=\emptyset, A \cap C \cap D=\emptyset$.
        \item $|F|\subseteq |A|\cap |B|$, $|F'|\subseteq |A|\cap |D|$,
        $|G|\subset |B|\cap |C|$, $|G'|\subset |D|\cap |C|$.
        \item The components of 
        $A=\sum n_j A_j$ and $C=\sum m_k C_k$ are disjoint and smooth, and the following equalities hold in $H^*(|A|)$ and $ H^*(|C|)$ respectively:
  \begin{align*}
        &[A]_{rel}\cdot [B]|_A=[F]_{rel},  &[A]_{rel}\cdot [D]|_{A}=[F']_{rel},  \\
        &[C]_{rel}\cdot [B]|_{C}=[G]_{rel},  &[C]_{rel}\cdot [D]|_{C}=[G']_{rel}. 
        \end{align*}
    \end{enumerate}
    Then, there is an equality:
    \[
        \int_A\langle B,C,D\rangle_{ABC}=\langle F,G'\rangle- \langle F',G\rangle.
    \]
\end{thm}
\begin{rem}
    Concretely, if $F=\sum m_k[F_k]$, the condition $(3)$ means $\sum_j n_j [B]|_{A_j}=\sum_{F_k\subseteq A_j} m_k[F_k]|_{A_j}$ for all $j$, and similarly for the other equalities.
\end{rem}
\begin{rem}
    One can likely drop the assumptions on the components of $A,C$ to be smooth and disjoint by working instead of $H^\ast(|A|)$, $H^\ast (|C|)$ with the analogous conditions in the relative cohomologies $H^\ast(Y,Y-|A|)$ and $H^\ast(Y,Y-|C|)$. For the purpose of this paper, the smooth version suffices.
\end{rem}
\begin{rem}
    As we realized during the final stages of preparing this article, the PhD thesis \cite[Thm 4.3.]{Schw_Diss} contains a result which is similar in spirit: It equates certain Massey products in Deligne cohomology (with integral coefficients), as defined by \cite{Den_higher}, with the difference of two height pairings. Unlike in our situation, the products considered in \cite{Schw_Diss} are always torsion.
\end{rem}
\begin{proof}
    Note that the integral vanishes on the indeterminacy ideal of the ABC Massey product because $F,F'$ are null-homologous on $Y$. In addition, the first hypothesis ensures that $F\cap G'=F'\cap G=\emptyset$. 
   
       Let $\omega_A$,...,$\omega_D$ be closed forms of pure type representing $[A]$,...,$[D]$, which are supported in small neighborhoods  of them chosen so that the triple intersections of those corresponding to $|A|,|B|,|C|$ and $|A|,|C|,|D|$ are empty. Further, choose $P_{BC}$ and $P_{CD}$ smooth forms of pure type such that $i\del\delbar P_{BC}=\omega_B\omega_C$ and $i\del\delbar P_{CD}=\omega_C\omega_D$. Moreover, choose currents $T$, $T'$ on $Y$, with $T|_{Y-|F|}=\dY\wedge \eta$, $T'|_{Y-|F'|}=\dY\wedge \eta' $ for some smooth forms $\eta,\eta'$ on $Y-|F|,Y-|F'|$ respectively, and such that $i\del\delbar T=\delta_F $, $i\del\delbar T'=\delta_{F'} $. Then, the left hand side is
    \begin{align*}
\int_A\omega_BP_{CD}-P_{BC}\omega_D=&\int_F P_{CD}- \int_{F'} P_{BC}= i\del \delbar T (P_{CD})-i\del \delbar T' (P_{BC})\\=&
\int_{Y} {\eta \omega_C \omega_D- \eta' \omega_B \omega_C}.
    \end{align*}
    In the first step, we used that $P_{CD}|_{A},P_{BC}|_A$ are $\del\delbar$-closed forms, due to the first hypothesis, and the first part of the last hypothesis. In the last, we used  $\supp(\omega_C) \subset (Y-|A|\cap|B|)\cap (Y-|A|\cap |D|)\subset (Y-|F|)\cap (Y-|F'|)$. Indeed, $\eta \omega_C$, $\eta' \omega_{C}$, are well defined $\del \delbar$-closed forms on $Y$. Using that
     $\eta|_{|C|}$, $\eta'|_{|C|}$ are also $\del \delbar$-closed smooth forms, and the second part of the last hypothesis, we obtain:
\[
\int_{Y} {\eta \omega_C \omega_D- \eta' \omega_B\omega_C} = \int_{G'} \eta -\int_{G}\eta'= \langle F,G' \rangle-\langle F',G \rangle.
\] 
\end{proof}

\begin{cor} \label{coro:blowup}
    Let $n=\dim(Y)=3$. Let $C_0,...,C_k$  be disjoint curves on $Y$ for some $k\geq 3$. Assume that $C_0,C_1,C_2,C_3$ have trivial normal bundle, $[{C_0}]=[{C_1}]$ and $[{C_2}]=[{C_3}]$. Let $\pi \colon \widetilde{Y}\to Y$ be the blow-up of $Y$ along $C_0\cup \ldots\cup C_k$ with exceptional divisor $E_i$ lying over $C_i$. Then, there is an equality
    \[
       \int_{E_0-E_1} \langle E_0+E_1,E_2+E_3,E_2-E_3\rangle = \langle C_0-C_1,C_2-C_3\rangle.
    \]
\end{cor}
\begin{proof}
 Since the normal bundle of $C_j$ are trivial for $j=0,1,2,3$, by choosing a non-vanishing section we can find a submanifolds $\bar{C}_j \subset E_j$ biholomorphic to $C_j$ via $\pi$. In addition, the equality $[E_j]^2=-[\bar{C}_j]=-\pi^*[C_j]$ makes the  ABC Massey product well defined. The hypotheses of Theorem \ref{thm: MP vs LN} hold by taking $F=\bar{C}_0-\bar{C}_1, G'=\bar{C}_2-\bar{C}_3, F'=G=0$. This yields
  $$
  \int_{E_0-E_1} \langle E_0+E_1,E_2+E_3,E_2-E_3\rangle = \langle \bar{C}_0-\bar{C}_1,\bar{C}_2-\bar{C}_3\rangle.$$
 The conclusion follows by observing that $\pi$ factors through the blow-up $\pi' \colon Y'\to Y$ of $Y$ along $C_0\sqcup C_1$ and applying  Proposition \ref{prop:adjunction} twice.
\end{proof}

\section{Construction of the manifold}
\subsection{Overview} We build the manifold $X=X_\tau$ as a resolution of singularities of an orbifold $\overline{X_\tau}$. The orbifold, in turn is built as the quotient of the torus \[
T:=T_\tau=\C^3/((\Z+\Z2\tau)\times\Z[i]^2)
\]
by an action of $(\Z/2\Z)^2=\langle \iota,\kappa\rangle$, where $\iota$ acts freely and $\kappa$ with fixed points. The intermediate space $M=T/\langle\iota\rangle$ is a mapping torus over $\C/(\Z+\Z\tau)$, illustrated in \Cref{fig: M}. The Massey product is built in the last subsection. It takes as inputs certain linear combinations of exceptional divisors on $X_\tau$, which lie over singularities arising from $\kappa$-invariant curves $N_0,\dots,N_3$ on $M$.

\subsection{The orbifold} Denote $\Lambda_\tau=\Z+\Z\tau\subseteq\C$ and $C_\tau=\C/\Lambda_\tau$. So $T=C_{2\tau}\times C_i\times C_i$. 

Define the involutions $\iota,\kappa:T\to T$ as follows:
\begin{align*}
    \iota([t,z_1,z_2])&=[t+\tau,-z_1,-z_2]\\
    \kappa([t,z_1,z_2])&=[-t,z_1,-z_2]
\end{align*}
and let $G=\langle \iota,\kappa\rangle\cong(\Z/2\Z)^2$. We want to compute the singular points of $T/G$, which amounts to a calculation of the locus on $T$ where $G$ acts with nontrivial stabilizer.

Because of the translation in the first coordinate, $\iota$ acts freely on $T$. Hence, the quotient $M:=T/\langle\iota\rangle$ is again a complex manifold. Projection induces a biholomorphism
\begin{equation*}
C_{2\tau}/(t\mapsto t+\tau)\cong C_\tau
\end{equation*}
and so the projection of $T$ to the first factor induces a holomorphic fibre bundle structure 
\begin{equation*}
    M\to C_\tau
\end{equation*}
with fibre the elliptic curve $C_i^2$. This projection is equivariant with respect to the induced action of $\kappa$ on $M$ and $(-1)$ on $C_\tau$. Hence, any fixed point of $\kappa$ on $M$ must lie over one of the four fixed points of $(-1)$. These are given by the set $A_\tau:=\{a_0^\tau,a_1^\tau,a_2^\tau,a_3^\tau\}$ of half-lattice points $a_0^\tau=[0]$, $a_1^\tau=[\frac{1}{2}]$, $a_2^\tau=[\frac{\tau}{2}]$, $a_3^\tau=[\frac{1+\tau}{2}]$. On the fibres over the fixed points $a_0^\tau,a_1^ \tau$, the induced action of $\kappa$ is given by $[z_1,z_2]\mapsto [z_1,-z_2]$. It has fixed points $C_i\times A_i$, where $A_i=\{a_0^i, a_1^ i,a_2^ i,a_3^i\}=\{[0],[\frac{1}{2}],[\frac{i}{2}],[\frac{1+i}{2}]\}$. On the other hand, since $[\tfrac{\tau}{2}]\overset{(-1)}\sim[-\tfrac{\tau}{2}]\overset{+\tau}\sim[\tfrac{\tau}{2}]$, on the fibres over the fixed points $a_2^\tau,a_3^\tau$, the induced action is given by $\kappa\circ\iota|_{C_i^2}:[z_1,z_2]\mapsto [-z_1,z_2]$, which has fixed points $A_i\times C_i$.

Let $\tilde{N}_{k,l}\subseteq T$ be the the following sets:
\begin{align*}
     \tilde{N}_{0,l}&=\left\{\left[0\right]\right\}\times C_i\times \{a_l^i\}\\
     \tilde{N}_{1,l}&=\left\{\left[\tfrac12\right]\right\}\times C_i\times \{a_l^i\}\\
     \tilde{N}_{2,l}&=\left\{\left[\tfrac\tau 2\right]\right\}\times \{a_l^i\}\times C_i\\
     \tilde N_{3,l}&=\left\{\left[\tfrac{1+\tau}2\right]\right\}\times \{a_l^i\}\times C_i
\end{align*}
and denote by $N_{k,l}:=q \tilde{N}_{k,l}$ the image under the projection $q:T\to M$. Then we can summarize the situation as follows:

\begin{lem}\,
\begin{enumerate}
    \item The points with nontrivial stabilizer of the $G$-action on $T$ have stabilizer equal to $\langle\kappa\rangle$. They consist of the disjoint union of the $\tilde{N}_{k,l}\cong C_i$ and their $\iota$-translates.
    
    \item       The fixed points of the $\kappa$-action on $M$ consist of the disjoint union of the $16$ elliptic curves $N_{k,l}\cong C_i$. Each of them has trivial normal bundle.

    \item Each $N_{k.l}$ (resp. $\tilde N_{k,l}$, resp. $\iota\tilde N_{k,l}$) has, as a fibre of the projection to $C_\tau$, trivial normal bundle in $M$ (resp. $T$).
\end{enumerate}
\end{lem}

\begin{figure}[h]
    \centering
    \includegraphics[width=0.7\linewidth]{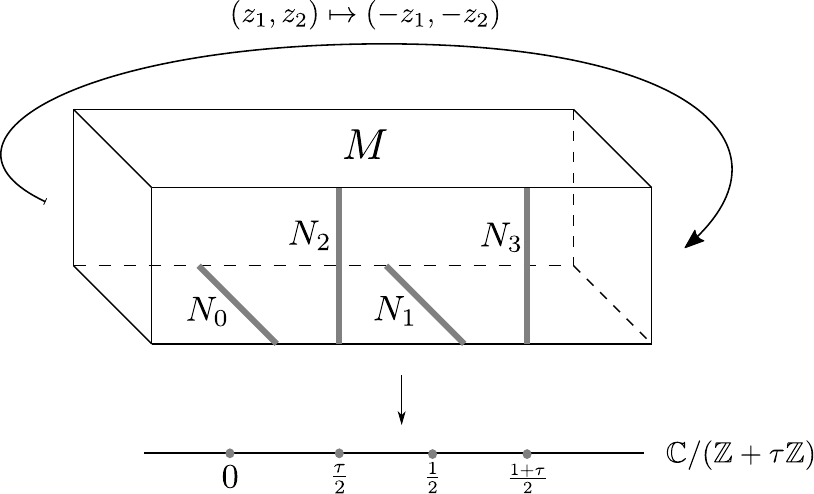}
    \caption{The mapping torus $M$ and the four curves $N_0,\dots,N_3$.}
    \label{fig: M}
\end{figure}

For later use, we abbreviate $N_{k,0}=:N_{k}$. The whole situation is illustrated in \Cref{fig: M}.

\subsection{The resolution}

We now describe a crepant resolution of the global quotient-type orbifold $\overline{X}:=\overline{X}_\tau:=T/G\cong M/\langle\kappa\rangle$ with trivial canonical bundle. First, consider the $\kappa$-equivariant map $\tilde{p}:\tilde{X}\to M$ given by blowing up all the $N_{k,l}$. The action of $\kappa$ on $\tilde{X}$ has the collection of exceptional divisors $E_{k,l}=\tilde{p}^{-1}(N_{k,l})$ as fixed points. To see this, note that the action on the normal bundle of each $N_{k,l}$ is given by multiplication by $(-1)$, which induces the trivial action on the projectivization. Since the $E_{k,l}$ are codimension one, the complex space underlying the orbifold $X:=\tilde{X}/\langle\kappa\rangle$ is a complex manifold and we obtain a resolution $p:X\to \overline{X}$.  This resolution is indeed crepant: Denoting $E=\sum E_{k.l}$, we have by standard formulas for the canonical bundle of blow-ups $K_{\tilde{X}}=\cO(E)$. On the other hand, $p$ is a ramified cover with ramification divisor equal to $E$, and so $p^*K_{X}\otimes \cO(E)=K_{\tilde{X}}$, which implies $K_X\cong\cO$.

\subsection{The cohomology of $X$}

By the blow-up formula for the cohomology of the blow up of compact K\"ahler manifolds in smooth centers (e.g. \cite[Thm 7.31]{Voisin}), the cohomology of $\tilde{X}$ can be computed as follows: For $h_{k,l}=c_1(\cO_{\mathbb{P}(\cN_{N_{k,l}})}(1))$ we have
\begin{equation}\label{eqn: HtildeX}
H(\tilde{X})\cong H(M)\oplus \bigoplus_{0\leq k,l\leq 3} H(N_{k,l}) h_{k,l}
\end{equation}
where the map from left to right is given by $\tilde p_*$ in the first component and the composition  \[H(\tilde{X})\to H(E)\cong H(\bigsqcup\mathbb{P}(\cN_{N_{k,l}}))\cong \bigoplus \left(H(N_{k,l})\oplus H(N_{k,l})h_{k,l}\right)\to H(N_{k,l})h_{k,l}\]
of restriction, projective bundle formula and projection.
The inverse is given by
\[
(a,bh)\mapsto \tilde{p}^*a + j_*(\tilde p|_E)^*b,
\]
where $j:E\to \tilde{X}$ is the inclusion.
Then $H(X)$ is calculated taking $\kappa$-fixed points in \eqref{eqn: HtildeX} and we obtain:
\begin{equation}\label{eqn: HX}
H(X)\cong H(\overline{X})\oplus \bigoplus_{0\leq k,l\leq 3} H(N_{k,l}) h_{k,l}
\end{equation}
with $H(\overline X)\cong H(T)^G\cong H(M)^{\kappa=\id}\cong\Lambda^r([dt],[d\bar{t}],[dz_1],[d\bar{z}_1],[dz_2],[d\bar{z}_2])^G$, which can, with its type-decomposition, be computed as follows
\[
H^{p,q}(\overline{X})\cong\begin{cases} \C &(p,q)=(0,0)\\
\langle [dtd\bar{t}], [dz_1d\bar{z}_1],[dz_2d\bar{z}_2]\rangle&(p,q)=(1,1)\\
\langle [dtdz_1dz_2]\rangle & (p,q)=(3,0)\\
\langle [dtdz_1 d\bar{z}_2],[dtdz_2d\bar{z_1}],[dz_1dz_2d\bar{t}]\rangle&(p,q)=(2,1)\\
\langle [d\bar td\bar z_1 d z_2],[d\bar td\bar z_2dz_1],[d\bar z_1d\bar z_2dt]\rangle&(p,q)=(1,2)\\
\langle [d\bar td\bar z_1d\bar z_2\rangle] & (p,q)=(3,0)\\
\langle [dtdz_1d\bar t d\bar z_1], [dtdz_2d\bar{t}d\bar z_2], [dz_1dz_2d\bar{z}_1d\bar z_2]\rangle& (p,q)=(2,2)\\
\langle [dtdz_1dz_2 d\bar t d\bar z_1d\bar z_2]\rangle&(p,q)=(3,3)\\
0&\text{else.}\end{cases}
\]

In particular, we get the following Hodge diamond for $X$:
\[
\begin{array}{ccccccccc}
& & & &1& & & &\\[3mm]
& & &0& &0& & &\\[3mm]
& &0& &19& &0& &\\[3mm]
&1& &19& &19& &1&\\[3mm]
& &0& &19& &0& &\\[3mm]
& & &0& &0& & &\\[3mm]
& & & &1& & & &\\[3mm]
\end{array}
\]

\subsection{Fundamental group}

We now prove that $X$ is simply connected. First, note that the map $p_*:\pi_1(X)\to \pi_1(\overline{X})$ is an isomorphism. This holds for resolutions in great generality \cite[Thm 7.5.]{Kollar}, \cite{tak} but can in this case also be seen via a Seifert-van-Kampen argument.

Denoting by $\bq:M\to\oX$ the quotient map, the induced map $\bq_*:\pi_1(M)\to \pi_1(\oX)$ is surjective because the action of $\kappa$ has a fixed point (see \cite[Chapter 2, Corollary 6.3]{Bre}). 
Thus, the simple connectedness of $X$ will follow once we prove that $\bq_*$ is trivial. 
To construct a set of generators of $\pi_1(M)$ we consider the short exact sequence
$$
\{1\} \to \pi_1(T,[0]) \to \pi_1(M,[0]) \to \{\iota, \id \} \to \{1\},
$$
where the first map is induced by the projection  $q \colon T \to M$, and the second  sends a loop in $M$ to the deck transformation that takes $[0]$ to the endpoint of its lift to $T$ starting at $[0]$.  Let $\gamma_1',\dots,\gamma_6' \colon [0,1]\to T$ be the canonical representatives of the generators of $\pi_1(T,[0])$. That is, $\gamma_1'(s)=[s,0,0]$, $\gamma_2'(s)=[2s\tau,0,0]$ and so on. Define 
$\gamma_k= q\circ \gamma_k'$ if $k\neq 2$ and $\gamma_2= q\circ \hat{\gamma}_2$, where $\hat{\gamma}_2\colon [0,1]\to T$, $\hat{\gamma}_2(s)=[s\tau,0,0]$. Then, $[\gamma_1],\dots, [\gamma_6]$ generate $\pi_1(M)$.

We finally show that $\bq \colon \pi_1(M) \to \pi_1(\oX) $ is trivial. We use repeatedly the following observation: if a loop $\gamma \colon [0,1] \to M$ satisfies $\kappa \circ \gamma= \gamma^{-1}$, then $\bq \circ \gamma$ is basepoint homotopic to a constant loop. This holds because in that case, $\bq \circ \gamma= \bq \circ \rho $, where 
$$
\rho(s)=
\begin{cases}
\gamma(s), &\mbox{if } t \leq 1/2, \\
\gamma^{-1}(s), &\mbox{if } t \geq 1/2,
\end{cases}
$$
and $\rho$ is homotopic to the constant loop via a basepoint-preserving homotopy in $M$.  Using this, together with the identities
$\kappa \circ \gamma_k'=\gamma_k'^{-1}$ for $k=1,5,6$ and $ \kappa \circ \iota \circ  \hat{\gamma}_2=\hat{\gamma}_2^{-1}$, we obtain $q_*[\gamma_k]=1$ if $k=1,2,5,6$.
To deal with the cases $k=3,4$, we consider the paths
 $$
 \tilde{\gamma}_2, \tilde{\gamma}_3, \tilde{\gamma}_4 \colon [0,1] \to T, \qquad \tilde{\gamma}_2(s)=[s\tau/2,0,0], \quad \tilde{\gamma}_3(s)=[\tau/2,s,0], \quad \tilde{\gamma}_4(s)=[\tau/2,is,0].
 $$
Observe that 
$
[\gamma_k']=[ \tilde{\gamma}_2 \cdot  \tilde{\gamma}_k \cdot  \tilde{\gamma}_2^{-1}]
$, 
so,  
$
\bq_*[\gamma_k]=(\bq\circ q)_* [\tilde{\gamma}_2 \cdot  \tilde{\gamma}_k \cdot  \tilde{\gamma}_2^{-1}]
$
. Since $ \kappa \circ \iota \circ \tilde{\gamma}_k=(\tilde{\gamma}_k)^{-1}$ it follows that
$
\kappa \circ (q \circ \tilde{\gamma}_k)=(q\circ \tilde{\gamma}_k)^{-1}
$, and thus
  $\bq \circ q \circ  \tilde{\gamma}_k$ 
  is base-point homotopic to a trivial loop. This implies 
$
(\bq\circ q)_* [\tilde{\gamma}_2 \cdot  \tilde{\gamma}_k \cdot  \tilde{\gamma}_2^{-1}]=1
$.

\subsection{Definition and evaluation of the Massey product}\label{sec: defn and evaluation}
 
  We recall there is a diagram
\[
X=\tilde{X}/\langle \kappa\rangle\overset{\pi}{\longleftarrow}\tilde{X}\longrightarrow M
\]
where the left hand map $\pi$ is the projection and the right hand map is the blow-up in the $N_{k,l}$. We recall $N_k=N_{k,0}\subseteq M$ and we write $E_k:=E_{k,0}$ for the exceptional divisor over $N_k$. Denote by $\bar{E}_k=\pi{E}_k\subseteq X$ the image of the exceptional divisors.  Let $D_1=\bar E_0-\bar E_1$, $D_2=\bar E_0 + \bar E_1$, $D_3:= \bar{E}_2+\bar E_3$, $D_4:= \bar{E}_2-\bar E_3$. 

With this notation, the expression referred to in \Cref{thmintro: ABC MP on CY} is
\[
\int_{D_1}\langle D_2, D_3, D_4\rangle_{ABC} .
\]

The proposition shows that to compute it, we can work on $\widetilde{X}$ instead of $X$.

\begin{prop} \label{prop:tildeX}
   
    There is an equality
    \[
   \int_{D_1}\langle D_2, D_3, D_4\rangle_{ABC} =8\int_{E_0-E_1}\langle E_0+E_1, E_2+E_3, E_2-E_3\rangle_{ABC}
    \]
\end{prop}
\begin{proof}
    We have $\pi_*[E_l]=[\bar E_l]$ and $\pi^*[\bar E_l]=2[E_l]$. Therefore:

\begin{align*}
\int_{D_1}\langle D_2, D_3, D_4\rangle_{ABC} &= \int_X[\bar E_0-\bar E_1]\cdot \langle \bar E_0+\bar E_1, \bar E_2+\bar E_3, \bar E_2-\bar E_3\rangle_{ABC}\\
&=\frac{1}{2}\cdot \int_{\widetilde{X}}\pi^*([\bar E_0-\bar E_1]\cdot \langle \bar E_0+\bar E_1, \bar E_2+\bar E_3, \bar E_2-\bar E_3\rangle_{ABC})\\
&=8\cdot \int_{\widetilde{X}}[{E}_0-{E}_1]\cdot \langle   E_0+ E_1,  E_2+ E_3,  E_2- E_3\rangle_{ABC}
\end{align*}
\end{proof}  

\begin{prop}[=\Cref{thmintro: ABC MP on CY}]\label{prop:lk-number-computations}
Let $\lambda\colon \C \to \C$ be the modular lambda function. There are equalities
$$
\int_{D_1}\langle D_2, D_3, D_4\rangle_{ABC}=8\cdot \langle N_{0}-N_{1}, N_{2}-N_{3}\rangle=\frac{4}{\pi} \log(|\lambda(\tau^{-1})|)= \frac{4}{\pi} \log(|1-\lambda(\tau)|).
$$
\end{prop}   
\begin{proof}
First, combining  \Cref{coro:blowup}  with  \Cref{prop:tildeX}, we find
\[
\int_{D_1}\langle D_2,D_3,D_4\rangle_{ABC}=8\cdot \langle N_0-N_1, N_2-N_3\rangle,
\]
and it remains to calculate the right hand side. We do so using the functoriality properties of the linking number to reduce to the case of a curve.

Consider $D=\{ [z_2]=[0]\}\subset M$ and denote by $j\colon D\hookrightarrow  M$ the inclusion. Since $N_1,N_2\subset D$ and the intersection of $D$ with $N_2,N_3$ is transverse and $j^*(N_2-N_3)=[\frac{\tau}{2},0,0]-[\frac{1+\tau}{2},0,0]$ we have by the second part of \Cref{prop:adjunction}:
$$\langle N_0-N_1, \left[\tfrac{\tau}{2},0,0\right]-\left[\tfrac{1+\tau}{2},0,0\right]\rangle_{D} = \langle N_0-N_1, N_2-N_3\rangle_M.$$
Considering now the identification $C_\tau\cong \{[z_1,z_2]=[0,0]\}\subset D, [x]\mapsto[x,0,0]$, a similar argument shows: 
\[\langle [0]-\left[\tfrac{1}{2}\right], \left[\tfrac{\tau}{2}\right]-\left[\tfrac{1+\tau}{2}\right]\rangle_{C_\tau}=\langle N_0-N_1, \left[\tfrac{\tau}{2},0,0\right]-\left[\tfrac{1+\tau}{2},0,0\right]\rangle_{D}.\]

Let now $\wp_\tau:C_\tau\to \CP^1$ denote the Weierstraß $\wp$-function for the lattice $\Z+\Z_\tau$. Since $\wp_\tau-\wp_\tau(\frac{1}{2})$  has $2[\tfrac{1}{2}]- 2[0]$ as divisor of zeros and poles,
example \ref{ex:linking-Riemann-Surface} shows:
\begin{align*}
\langle  [0]-\left[\tfrac{1}{2}\right],\left[\tfrac{\tau}{2}\right]-\left[\tfrac{1+\tau}{2}\right] \rangle_{C_\tau} &= 
 \frac{1}{2\pi } \log\left( \dfrac{|\wp_\tau (\frac{1+\tau}{2})- \wp_\tau (\frac{1}{2})|}{|\wp_\tau (\frac{\tau}{2})- \wp_\tau (\frac{1}{2})|}\right).
\end{align*}
The statement now follows from the definition and properties of the modular lambda function \cite[Chapter 7.3]{Ahlfors}, namely: 
$$
\dfrac{\wp_\tau (\frac{1+\tau}{2})- \wp_\tau (\frac{1}{2})}{\wp_\tau (\frac{\tau}{2})- \wp_\tau(\frac{1}{2})}= \lambda(\tau^{-1})= 1-\lambda(\tau).
$$
\end{proof}

\begin{rem} Let us outline two different ways of doing the computation of the linking number $\langle N_0-N_1, N_2-N_3\rangle$:
\begin{enumerate}
    \item Instead of using using the functoriality for embeddings, we may use it for submersions (so the first part of \Cref{prop:adjunction} to again reduce to a curve:
    Note that the fibration $M \to C_\tau$ factors as a composition of two holomorphic fibrations
$M \to M_2 \to C_\tau$, where $M_2=(C_{2\tau}\times C_i)/\iota_2$ with $\iota_2[t,z]=[t+\tau,-z]$. The projections involved are $p_1[t,z_1,z_2]_{\iota}=[t,z_1]_{\iota_2}$ and $p_2[t,z]_{\iota_2}=[t]$. For $k=0,1$ set $N_{k}^2=\{ [\frac{k}{2},z]_{\iota_2} \}$.

Observe that $p_1\colon N_{0}\cup N_{1} \to N_{0}^2\cup N_{1}^2$ is a diffeomorphism and $p^{*}_1([\frac{\tau}{2},0]-[\frac{1+\tau}{2}])=N_{2}-N_{3}$. In addition, 
$p_2$ restricts to a diffeomorphism from $\{[\frac{\tau}{2},0]\}\cup\{[\frac{1+\tau}{2},0]\}$ to $\{\frac{\tau}{2}, \frac{1+\tau}{2} \}$ and $p_2^*([0]-[\frac{1}{2}])= N_{0}^2-N_{1}^2 $. The dimensional hypotheses in Proposition \ref{prop:adjunction} hold; applying the first part twice, we obtain:
$$
\langle N_{0}-N_{1}, N_{2,}-N_{3}\rangle_M =\langle  [0]-\left[\tfrac{1}{2}\right],\left[\tfrac{\tau}{2}\right]-\left[\tfrac{1+\tau}{2}\right] \rangle_{C_\tau}
$$
and one proceeds as before.
\item One can also do the computation on $M$ directly: Let $f:D\to \CP^1$ be the composition of inclusion, projection and the shifted Weierstraß $\wp$-function \[D\overset{j}{\hookrightarrow}  M\longrightarrow C_\tau\overset{\wp_\tau-\wp_\tau(\frac{1}2)}{\longrightarrow}\mathbb{CP}^1.\] Then, $T:=j_*(\frac{\log \|f\|}{2\pi})$ is a current which satisfies $i\del\delbar T= \delta_{N_0}-\delta_{N_1}$. Using \cite[Ex. 8.2.5]{Hoermander_LPDO_I}, one may check that it satisfies the conditions of \Cref{lem:flexibility} and one can then evaluate 
\[\langle N_0-N_1,N_2-N_3\rangle=(\delta_{N_2-N_3}\wedge T)(1)=\frac{1}{2\pi}j_*(\log|f|\cdot \delta_{N_2\cap D-N_3\cap D)})(1)
\]
using that $N_i\cap D$ consists of a single point for $i=2,3$.
\end{enumerate}
\end{rem}
\begin{rem}
Similar computations shows that for every pair $0\leq l,l'\leq 3$, we have
$$
\langle N_{0,l}-N_{1,l}, N_{2,l'}-N_{3,l'}\rangle=\frac{1}{2\pi} \log(|\lambda(\tau^{-1})|)= \frac{1}{2\pi} \log(|1-\lambda(\tau)|),
$$
yielding to more non-vanishing ABC Massey products.
\end{rem}
\begin{rem}
In the computation of the ABC Massey product in the examples in \cite[Theorem 5]{PSZ24}, we worked with forms that approximate some naturally occuring currents (see discussion directly after the statement of \cite[Theorem 5]{PSZ24}).
By using \Cref{thm: MP vs LN} and \Cref{lem:flexibility} of the present paper, one can avoid this approximation and directly argue with the natural currents, as in the preceding computations.
\end{rem}
\bibliographystyle{acm}
\bibliography{biblio}

\end{document}